\documentclass[twoside,reqno]{lt6-proc}
\usepackage{epsfig,cite}
\usepackage{amssymb,amsmath}
\usepackage{times}
\begin{document}
\sloppy \raggedbottom
\setcounter{page}{1}

\newpage
\setcounter{figure}{0}
\setcounter{equation}{0}
\setcounter{footnote}{0}
\setcounter{table}{0}
\setcounter{section}{0}

\title{What can one reconstruct from the representation ring of a 
compact group?}

\runningheads{Zolt\'an Zimbor\'as}{Reconstruction from Representation Rings}

\begin{start}

\author{Zolt\'an Zimbor\'as}{1}

\address{Research Institute for Particle and Nuclear Physics, Budapest \\ 
H-1525, P.O. Box 49, Hungary, e-mail:cimbi@rmki.kfki.hu}{1}

\begin{Abstract}
It is well known that there exist non-isomorphic compact groups
with isomorphic representation rings (fusion rules).
Nevertheless, considerable structural information about the group
can be reconstructed from its representation ring.
We review these types of partial reconstruction theorems, including
some recent results. In the Appendix a 
derivation of the Clebsch-Gordan series
of $SU(2)$ based only on information about the dimensions of the irreducible
representations is presented.
\end{Abstract}
\end{start}


\section{Introduction}

Although the representation theory of compact groups is
far from being completely understood, 
it is a very well established
field of mathematics. 
One of the basic theorems in this field is that every continuous 
unitary representation of a compact group on a complex Hilbert space
is equivalent to a 
direct sum of irreducible finite-dimensional
unitary representations (shortly: irreps). 
In particular, the tensor product of any
two irreps can be decomposed to the so called {\it Clebsch-Gordan series}: 
\begin{equation*}
D_{p} \otimes D_{q} \cong  \bigoplus_{r}  N^{r}_{p,q} D_{r}. 
\end{equation*} 
The multiplicities $N^{r}_{p,q}$ are often referred to as 
{\it fusion rules}.
In terms of irreducible characters the above decomposition can 
be written in the form   
\begin{equation}
\chi_{p} \cdot \chi_{q}= \sum_{r} N^{r}_{p,q}\chi_{r}, \quad \chi_{r} 
\in I_{G}, \label{CGchar}
\end{equation}
here and throughout the paper $I_{G}$ will denote the set of 
characters 
belonging to the irreps of a compact group $G$.

From equation \eqref{CGchar} the pointwise product of 
characters of arbitrary  
finite-dimensional continuous unitary
representations can be deduced (provided that their  
decompostions to irreducibles are known),
and one can also extend this product to {\it generalized characters}, i.e. 
to differences of characters.
In this way a ring $\mathcal{R}(G)$ is obtained. This ring
admits a natural partial ordering making it into a partially ordered
ring with positive cone $\mathcal{R}(G)^{+}$ consisting of 
characters, i.e., $\chi_{2} \prec \chi_{1}$ holds, if 
$(\chi_{1}-\chi_{2}) \in \mathcal{R}(G)^{+}$. 
The ordered ring $\mathcal{R}(G)$ is called the {\it representation ring}
of the compact group $G$. The irreducible characters are the minimal
positive elements of $\mathcal{R}(G)$ 
and form a $\mathbb{Z}$-basis for the ring.
The representation ring 
encodes exactly the same 
information as the fusion rules.

Considerable amount of work has been devoted to
the derivation of fusion rules of various compact groups.  
In this contribution we will consider the
less frequently investigated {\it dual question}:
Provided that the fusion rules, or representation ring,
of a (possibly unknown) compact group
is known, what attributes of the
group can one reconstuct from this information?
Our intention is 
to give a 
short and informal review on this subject.

\section{Complete reconstructions}

In general it is
not possible to recover a compact group completely
from its representation ring, 
since there exist non-isomorphic compact groups
- even finite groups - with isomorphic representation rings (see Fig. 1.).
\begin{figure}[h]
\centerline{\epsfig{file=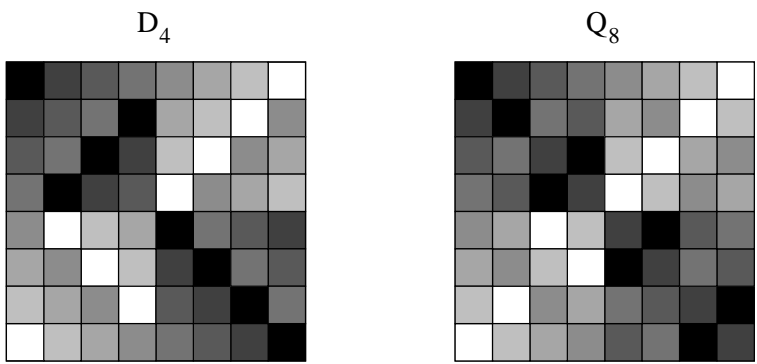 , width=100mm}}
\caption{Multiplication tables of 
the dihedral group $D_{4}$ and the quaternion group $Q_{8}$,
in both tables black squares denote the unit elements.
These two groups 
are not isomorphic, 
since the number of elements of order two differs,
but the groups have equivalent character tables, thus 
isomorphic fusion rules.}\label{D4Q8}
\end{figure}

It is quite remarkable, 
that the phenomenon illustrated in Fig. 1.
cannot occur in the case of {\it connected} compact groups. 
According to a beautiful theorem of Handelman \cite{Handelman},
two connected compact groups can only have
order isomorphic representation rings (equivalent fusion rules),
if the groups themselves are isomorphic.
However, Handelman's theorem is not a reconstruction theorem -
it only states this unicity,
but does not provide a general reconstruction method.
 
There are, however, other types of
representation theoretical data that can be used
to completely recover any compact group.
Namely, a compact group can be reconstructed from  
its Krein-algebra (Tannaka-Krein reconstruction theorem \cite[chapter 7]{HR}),
or from its representation category, 
viewed as a rigid, monoidal, symmetric
C$^{*}$-category (Doplicher-Roberts reconstruction theorem 
\cite{DR}).
This latter theorem 
provides also the basis for 
the theory of superselection sectors in Algebraic Quantum Field Theory 
\cite{Haag}.

\section{Pontryagin duality and the abelianization of a compact group}

Let $G$ be a locally compact abelian group, and let $\widehat{G}$
denote the set of continuous unitary one-dimensional characters.
The product of two such characters will again be in $\widehat{G}$, 
hence $\widehat{G}$ has a natural abelian group structure.
Furthermore, it can also naturally be endowed with a topology 
\cite[chapter 6]{HR}, which  
is compatible with its group structure.
This topology of $\widehat{G}$ 
is simply the discrete topology, when $G$ is compact. 
$\widehat{G}$, regarded as a topological group, is called the dual of $G$.
Also the dual of $\widehat{G}$ can be constructed,
the following well-known theorem shows the importance of this
"double dual" construction:\\

{\raggedright
{\bf Theorem } (Pontryagin-van Kampen duality theorem)\\
{\it Let $G$ be a locally compact abelian group, and let 
$\widehat{\widehat{G}}$ be the dual of its dual. $\widehat{\widehat{G}}$
and $G$ are topologically isomorphic.}}
\\

This theorem gives a method to recover an 
abelian compact group $G$ from its representation ring.
The fusion rules in this case simply encode the multiplication laws
of $\widehat{G}$, and as  
mentioned earlier, $\widehat{G}$ has a trivial 
discrete topology.
Now, according to the Pontryagin-van Kampen theorem, 
$G$ can be recovered (up to topological isomorphism)
by finding the 
one-dimensional unitary characters of $\widehat{G}$.

If $G$ is a nonabelian compact group, then the group of 
continuous, unitary one-dimensional characters of $G$ will be
isomorphic to the dual of the {\it abelianization} of $G$, i.e. to
the dual of $G/[G,G]$, where $[G,G]$ 
denotes, as usual, the commutator subgroup of $G$ 
\cite[Theorem 23.8]{HR}. 
Hence, the
abelianization of $G$ can (in principle) be
recovered from $\mathcal{R}(G)$.

\section{The center of a compact group and the chain group construction}

The other abelian compact group, beside the
abelianization, that can canonically be
associated to every compact group $G$ is its {\it center}.
It is natural to ask whether, similarly to the abelianization, 
also the center can be reconstructed from $\mathcal{R}(G)$. 
This question 
has only recently been answered by
Baumg\"artel and Lled\'o \cite{BL}, and M\"uger \cite{M}.

In \cite{BL}
the so called chain group construction was introduced.
The construction 
is based on the following equivalence relation 
$\sim$ on the set of irreducible characters:   
let $\chi, \eta \in I_{G}$, 
we define $\chi \sim \eta$, if there exists 
$\psi_{1}, \psi_{2}, \dots , \psi_{n} \in I_{G}$ (a "chain" of
irreducible characters), such that 
$\chi, \eta \prec \psi_{1} \cdot \psi_{2} \cdot \ldots \cdot \psi_{n}$,
i.e., both $\chi$ and $\eta$ is contained 
in the irreducible decomposition of the above product. 

The $\sim$-equivalence class of a character $\chi$ will be denoted by 
$\langle \chi \rangle$. The quotient $\mathcal{C}(G) := I_{G}/\sim$
is an abelian group, with respect to the operations 
$\langle \chi \rangle \cdot \langle \psi \rangle = \langle \eta \rangle$,
where $\eta$ is any irreducible character with the property
$\eta \prec \chi \cdot \psi$. The group $\mathcal{C}(G)$ is called
the {\it chain group} of $G$.

As an example, let us derive the chain group of $SU(2)$.
Let $\chi_{1/2} \in I_{SU(2)}$ denote
the unique two-dimensional irreducible character of $SU(2)$.
For any two odd-dimensional irreducible characters $\psi , \eta \in I_{SU(2)}$,
there is an appropriately large even number $n$, such that 
$\psi, \eta \prec \chi_{1/2}^{n}$. Similarly, for any  
two even-dimensional irreducible characters, there is an approriately
large odd number $m$, such that both even-dimensional characters 
are contained in $\chi_{1/2}^{m}$. However, no products of
irreducible characters contains 
simultaneously an even- and an odd-dimensional element of $I_{SU(2)}$.
Thus, $SU(2)$ has exactly two $\sim$-equivalence classes,
and $\mathcal{C}(SU(2)) \cong Z_{2}$.
 
The meaning of the chain group is revealed by
the following theorem, which was conjectured by 
Baumg\"artel and Lled\'o,
and was proved by M\"uger \cite{M}:\\

{\raggedright {\bf Theorem}\\
{\it The chain group $\mathcal{C}(G)$ of a compact group $G$ 
and $\widehat{Z(G)}$, the 
dual of the center of $G$, are isomorphic.}
}
\\

This means that one can reconstruct the dual of the center
of $G$ from $\mathcal{R}(G)$, and thus in principle also
$Z(G)$.

This reconstruction theorem can also be used to
prove several propositions about the center of
various types of compact groups. 
An example of such a theorem is the following: 
if the non-zero fusion rules of a compact group
are odd numbers (e.g. if the group is simply reducible), 
then the existence of a pseudo-real
(symplectic) representation implies that the group has 
a non-trivial center.
This theorem is a consequence 
of an interesting fusion property of pseudo-real representations \cite{FRS}.
The detailed proof of this and similar propositions together with a
generlization of the chain group construction will be published in
\cite{ZZ}.

\section{Closed normal subgroups and representation subrings}
 
A subring $R$ of $\mathcal{R}(G)$ will be called a
{\it representation subring}, if it is spanned by irreducible 
characters, 
i.e., if $\exists$ $J \subset I_{G}$ which is a $\mathbb{Z}$-basis of $R$.
There is a one-to-one Galois-correspondence between  
representation subrings $R \subset \mathcal{R}(G)$
and closed normal subgroups $H \subset G$ given by
\begin{eqnarray*}
H_{R} &=& \{ g \in G \bigm|  g \in Ker(\chi), \; \forall \chi \in R^{+} \}, \\
R^{+}_{H} &=& \{ \chi \in \mathcal{R}(G)^{+} \bigm| 
g \in Ker(\chi), \; \forall g \in H  \},
\end{eqnarray*}
where $R^{+}$ denotes the positive elements of a representation subring
$R$, i.e., the characters contained in $R$. 
This means that one can reconstruct the number (or cardinality) 
of closed normal subgroups of a compact group from its representation ring.
Furthermore,
\begin{equation*} 
R_{H} \cong \mathcal{R}(G/H), 
\end{equation*}
hence also the representation rings of the quotient groups can be obtained.

Moreover, for every $R_{H}$ the chain group construction can be 
generalized in a non-trivial way \cite{ZZ}, and one 
can recover also the dual of $Z(G) \cap H$.
This means that the number 
(or cardinality) of closed normal subgroups 
that contain a given subgroup of the center can be 
reconstructed, too. 

It is noteworthy that the 
representation subring belonging to the center, i.e. $R_{Z(G)}$, can be
canonically characterized for every compact $G$, it is
the smallest representation subring that
contains all irreducible characters that
appear in a decomposition of $\chi \cdot \chi^{*}$ for some 
$\chi \in I_{G}$ 
\cite{ENO, M}.

For the sake of completeness let us also mention that if $G$ is 
a {\it connected} compact group, and $R_{H}$ is a 
representation subring of $\mathcal{R}(G)$, then
one can recover from $\mathcal{R}(G)$ whether
$H$ is invariant with respect to all automorphisms of $G$. 
Namely,
$H$ has this property
if and only if all order automorphism of $\mathcal{R}(G)$
leaves $R_{H}$ invariant
(this is a direct consequence of
Theorem 2.14. in \cite{Handelman}). 
For non-connected compact groups
this does not hold in general.

\section{Character tables of finite groups and the dimensions of 
irreps of compact groups}

The way to obtain fusion rules from
the character tables of finite groups is a subject which is 
covered by almost all standard textbooks on finite groups.
It is, however, seldom mentioned 
that also the character table of a finite group can be 
recovered from the
fusion rules. This can be done, because the only non-trivial 
sets of complex solutions $\{ \alpha_{p} \}$ of the set of equations
\begin{equation}
\alpha_{p} \cdot \alpha_{q} = \sum_{r} N_{p,q}^{r} \alpha_{r} 
\label{chartabl}
\end{equation}
are exactly the columns of the character table, i.e. the character values
of certain conjugacy classes. This is not true for
compact groups in general.

From the character table  
of a finite group one can directly obtain
the abelianization, the center, and the number of normal subgroups, 
i.e. properties that can also be reconstructed from
the fusion rules of a compact group in general. 
In addition, one can recover the dimensions of the irreps,
the number of conjugacy classes that are left invariant
by the inverse, the sizes of the conjugacy classes,
and "the conjugacy class structure" of the 
commutator subgroup and other normal subgroups, 
and much other similar information \cite{Lomont}.

In the case of generic compact groups the set of equations
equations \eqref{chartabl} has usually 
many other solutions than 
the character values of a given group element (or conjugacy class).
One can not even obtain in general the dimensions of the irreps
from \eqref{chartabl}, since there might be several sets of postive integers
that satisfy this equation (this is the case even for $SU(2)$).
Despite this, there 
is a more complicated way of obtaining the dimension of the irreps, 
which is presented in \cite{Handelman}. 
It is an open question, at least according to the
authors knowledge, whether also other character values
can be recovered from the representation ring for general
compact groups. It is interesting that in the case of $SU(2)$
one can recover the fusion rules only from the dimensions of the irreps,
we present this short and easy derivation in the Appendix.

\section{Appendix - Deriving the Clebsch-Gordan series of SU(2) 
from the dimensions of the irreducible representations}

There exist many 
methods for deriving the Clebsch-Gordan series
of $SU(2)$ (for a recent and very nice derivation see \cite{Jones}).
We present a new type of derivation, which is based only
on the compactness of $SU(2)$ and the property that it
has one and only one equivalence class of irreps 
for every dimension. 
It is noteworthy that the only compact group with 
this latter property is $SU(2)$.

As usually, $D_{j}$ ($j \in \frac {1}{2} \mathbb{Z}$)
will denote the $2j+1$ dimensional (up to unitary equivalence) unique
irreducible representation
of SU(2), and $\chi_{j}$ denotes its character.
Due to the "uniqueness property"
$D_{j}$ must be self-dual 
(equivalent to its dual)
, thus the decomposition  of $D_{j} \otimes D_{j'}$ contains 
the trivial representation $D_{0}$ only if $j=j'$, 
and even then only once, i.e. $N^{0}_{j,j}=1$.

First, the decompositions of
tensor products of type 
$D_{1/2} \otimes D_{j}$ will be derived. 
Consider the first non-trivial among these, 
$D_{1/2} \otimes D_{1/2}$. This is a 
4-dimensional representation and, as already mentioned,
it contains $D_{0}$ exactly once. This and the fact that
the dimensions should match allow only one possible decomposition:
$D_{0} \oplus D_{1}$.
The more general formula
\begin{equation} 
D_{1/2} \otimes D_{j} \cong D_{j-1/2} \oplus D_{j+1/2}, \label{CG1}
\end{equation}
can now be proved by induction.
Suppose that \eqref{CG1} holds for any non-zero $j$
less than a specific $k \in \frac {1} {2} \mathbb{Z}$. 
We will prove it also for $k$.
If some $D_{j}$ is contained in the decomposition of
$D_{1/2} \otimes D_{k}$
with a non-zero multiplicity, then $D_{0}$ has to be
conatined in 
$D_{j} \otimes D_{1/2} \otimes D_{k}$ with the same multiplicity.
But \eqref{CG1} was supposed to hold  for $0 < j < k $, thus
in this case we can rewrite this product as  
$(D_{j-1/2} \oplus D_{j+1/2}) \otimes D_{k}$.
This representation contains $D_{0}$
only when $j=k - \frac {1} {2}$ (for $j<k$), i.e.,
the only representation with $j<k$ that is contained in 
$D_{1/2} \otimes D_{k}$ is $D_{k-1/2}$.
Hence the only possible decomposition of $D_{1/2} \otimes D_{k}$,
such that the dimensions match and the aboved-mentioned
condition holds is $D_{k-1/2} \oplus D_{k+1/2}$, and \eqref{CG1} is proved by 
induction. This was the most important part of our derivation.
 
Now, as the decompositon of the tensor product of 
$D_{1/2}$ and any other irrep of $SU(2)$
is known, one can decompose any power of $D_{1/2}$:
\begin{equation*}
D_{1/2}^{\otimes n} = D_{n/2}\oplus 
\bigoplus_{i=1}^{[n/2-1]}\left[ {n \choose i} -{n \choose i-1} \right]
D_{n/2-i} ,
\end{equation*}
where $[n/2-1]$ denotes the integer part of $n/2-1$.
The proof is left to the reader. 
In terms of characters this reads as:
\begin{equation}
\chi_{1/2}^{n} = \chi_{n/2}+
\sum_{i=1}^{[n/2-1]}\left[ {n \choose i} -{n \choose i-1} \right]\chi_{n/2-i}. \label{CG2}
\end{equation}

From equation \eqref{CG2} one can express any irreducible characters in terms
of $\chi_{1/2}$ powers (note that $\chi_{1/2}^{0}=\chi_{0}$):
\begin{equation}
\chi_{j}= \chi^{2j}_{1/2} + \sum^{[j]}_{i=1}(-1)^{i}{2j-i \choose i }
\chi^{2j-2i}_{1/2}.  \label{CG3}
\end{equation} 
The product $\chi_{j}\cdot \chi_{j'}$ can now be evaluated
using equations \eqref{CG2} and \eqref{CG3}, the result is:
$\chi_{j} \cdot \chi_{j'}= \chi_{|j-j'|} + \chi_{|j-j'|+1} + \cdots+ 
\chi_{j+j'-1} + \chi_{j+j'}$.
Converting this back to representations one recovers the usual Clebsch-Gordan
series of SU(2):
\begin{equation*}
D_{j} \otimes D_{j'} \cong D_{|j-j'|} \oplus D_{|j-j'|+1} \oplus \cdots
\oplus D_{j+j'-1} \oplus D_{j+j'}.
\end{equation*}


\section*{Acknowledgments}

I would like to thank the Organizers for 
having arranged this nice Conference, and I 
would also like to thank
the financial support by EPS and the support for travel by EMS.
Useful discussions with Sz. Farkas, P. P. P\'alfy, and
P. Vecserny\'es are thankfully acknowledged.  
This work was partly supported by grant OTKA T043159.


\end{document}